\documentclass[12pt]{article}
\usepackage{amssymb}
\usepackage{latexsym}
\usepackage{epsfig}
\usepackage{graphicx}
\usepackage{psfrag}

\newcommand{\newtekst}{\bf}
\newcommand{\newtext}{}

\newtheorem{definition}{Definition}[section]
\newtheorem{theorem}[definition]{Theorem}

\newcommand{\R}{{\mathbb R}}

\renewcommand{\S}{{\mathbb S}}

\newcommand{\proofbox}{{$\Box$}}

\def\tilde{\widetilde}
\def \bfo {\begin {eqnarray*} }
\def \efo {\end {eqnarray*} }
\def \ba {\begin {eqnarray*} }
\def \ea {\end {eqnarray*} }
\def \beq {\begin {eqnarray}}
\def \eeq {\end {eqnarray}}
\def \supp {\hbox{supp}\,}

\def \dist {\hbox{dist}}
\def\diag{\hbox{diag }}
\def \det {\hbox{det}}

\def \e {\varepsilon}
\def \p {\partial}
\def \a {\alpha}

\def\p{\partial}
\def\sph{\mathbb S^2}
\def\R{\mathbb R}

\title{Electromagnetic wormholes via \\
       handlebody constructions}

\author{Allan Greenleaf, Yaroslav Kurylev,
\\ Matti Lassas and Gunther Uhlmann
\thanks{AG and GU are supported by
US NSF, ML by CoE-programm 213476 of Academy of Finland.}}

\date{}

\begin{document}

\maketitle

\begin{abstract}
Cloaking devices are prescriptions of electrostatic, optical or electromagnetic
parameter fields (conductivity $\sigma(x)$, index of
refraction $n(x)$,  or electric permittivity $\epsilon(x)$ and magnetic
permeability $\mu(x)$) which are
piecewise smooth on $\R^3$ and singular on a  hypersurface
$\Sigma$, and such that objects in
the region enclosed by $\Sigma$ are not detectable to external
observation by waves. Here, we give
related constructions of invisible tunnels, which allow
electromagnetic waves to
pass  between possibly distant points, but with only the ends of the
tunnels  visible
to electromagnetic imaging. Effectively, these change the topology of
space with
respect to solutions of Maxwell's equations, corresponding to attaching a
handlebody to $\R^3$. The resulting devices thus function as electromagnetic
wormholes.
\end{abstract}

\vfil\eject

\section{Introduction}

There has recently been considerable interest, both theoretical
\cite{Le,PSS1,PSS2, GKLU,KSVW}
and experimental \cite{SMJCPSS}, in invisibility (or ``cloaking") from
observation by electromagnetic (EM) waves. (See also \cite{MBW} for a
treatment of
cloaking in the context of elasticity.) Theoretically, cloaking
devices  are given by
specifying the conductivity $\sigma(x)$ (in the case of electrostatics), the
index of refraction
$n(x)$ (for optics in the absence of  polarization, where one uses
the Helmholtz equation), or the electric permittivity $\epsilon(x)$
and magnetic
permeability $\mu(x)$ (for the full system of Maxwell's equations.)
In the constructions
to date, the EM parameter fields ( $\sigma;  n;\epsilon\hbox{ and }\mu$ ) have
been piecewise smooth
and anisotropic. (See, however, \cite[~Sec.4]{GLU} for an example
that can be interpreted
as cloaking with respect to Helmholtz by an isotropic negative index
of refraction
material.) Furthermore, the EM parameters have singularities, with one or more
eigenvalues of the tensors going to zero or infinity as one approaches from on
or  both sides the {\it
cloaking surface} $\Sigma$, which encloses the region within which
objects may be  hidden from
external observation.
Such constructions might have remained theoretical curiosities, but
the advent of
{\it metamaterials}\cite{E} allows one, within the constraints of current
technology, to construct media with fairly arbitrary $\epsilon(x)$
and $\mu(x)$.

It thus becomes an interesting mathematical problem with practical
significance to
understand what other new phenomena of wave propagation can be
produced by prescribing
other arrangements of $\epsilon$ and $\mu$.
Geometrically,  cloaking
can be  viewed as arising from a singular transformation of $\R^3$.
Intuitively,  for a spherical cloak \cite{GLU1,GLU2,PSS1}, it is as
if an infinitesimally small hole in  space
has been stretched to a ball $D$; an
object can be inserted inside
the  hole so created and is then invisible to external observations.
On the level of the EM parameters, homogeneous, isotropic parameters
$\epsilon,\mu$ are pushed
forward to become inhomogeneous, anisotropic and singular as one
approaches $\Sigma=\p D$
from the exterior. There are then two ways,  referred to
as the {\it single} and {\it double} coating in \cite{GKLU}, of continuing
$\epsilon,\mu$ to within $D$
so as to rigorously obtain invisibility with respect to locally
finite energy waves.
We refer to either process as {\it blowing up a point}. As observed
in \cite{GKLU},
one can use the double coating to produce a manifold with a different topology,
but with the change in topology invisible to external measurements.

To define the solutions of Maxwell's equations rigorously
in the single coating case, one has to add boundary conditions on
$\Sigma$. Physically, this corresponds to the lining of the interior
of the single coating material, e.g., in the case of
blowing up a point, with a perfectly conducting layer,
see  \cite{GKLU}.
  We point out here that in the recent preprint \cite{Weder}, the
single coating construction is supplemented with
selfadjoint extensions
of Maxwell operators in the interior of the cloaked regions;
these implicitly impose interior boundary conditions on the boundary of the
cloaked region, similar
to the  PEC boundary condition suggested in \cite{GKLU}. For the case of an
infinite cylinder the Soft-and-Hard (SH) interior boundary condition is used in
\cite{GKLU} to guarantee cloaking of active objects, and is   needed even for
passive ones.

In this paper, we show  how
more elaborate geometric constructions, corresponding to {\it blowing up a
curve}, enable the description of tunnels which allow the passage of waves
between distant points, while only  the ends of the tunnels
are visible to external observation. These devices function as
electromagnetic wormholes, essentially changing the topology of
space with respect to
solutions of Maxwell's equations.

{\newtext
We form the wormhole device around an obstacle $K\subset\R^3$ as follows.
First, one surrounds $K$ with metamaterials, corresponding to a specification
of EM parameters $\tilde{\e}$ and $\tilde{\mu}$. Secondly, one lines the
surface of $K$ with material implementing the Soft-and-Hard (SH)  boundary
condition from antenna theory \cite{HLS, Ki, Ki2}; this condition arose
previously \cite{GKLU} in the context of cloaking an infinite cylinder.  The EM
parameters, which become singular as one approaches $K$, are given as the
pushforwards of nonsingular parameters $\e$ and $\mu$ on an abstract
three-manifold $M$, described in Sec. 2. For a curve $\gamma\subset M$, we
construct the  diffeomorphism $F$ from $M\setminus\gamma$ to the
wormhole device in Sec. 3. For the resulting EM parameters $\tilde{\e}$ and
$\tilde\mu$, we have singular coefficients of Maxwell's equations at $K$, and
so it is necessary to  formulate an appropriate notion of locally finite energy
solutions (see Def. 4.1). In Theorem 4.2, we then show that there is a perfect
correspondence between the external measurements of EM waves
propagating through
the wormhole device and those propagating on the wormhole manifold.}

It was shown in
\cite {GKLU} that the cloaking constructions are mathematically valid
at all frequencies
$k$. However, both cloaking and  the wormhole effect studied here
should be considered as
essentially monochromatic, or at least very narrow-band, using
current technology, since,
from a practical point of view  the metamaterials needed to implement
the constructions
have to be fabricated and assembled with a particular wavelength in mind, and
theoretically are subject to significant dispersion \cite{PSS1}.
Thus, as for cloaking in \cite{Le,PSS1,GKLU}, here we describe the
wormhole construction relative to electromagnetic waves at a fixed positive
frequency
$k$. We point out that the metamaterials used in the experimental verification
of cloaking \cite{SMJCPSS} should be readily adaptable to yield a physical
implementation, at microwave frequencies, of the wormhole device 
described here.
See Remark 1 in Sec. \ref{four-two} for further discussion.

The results proved here were announced in \cite{GKLU2}.

\section{The wormhole manifold $M$}

First we explain, somewhat informally,  what we mean by a  wormhole.
The concept of a wormhole  is familiar
from general relativity \cite{T1,T2}, but here we define a wormhole
as an object
obtained by  gluing together pieces of Euclidian space equipped with
certain anisotropic EM
parameter fields. We start by describing  this process heuristically;
later, we explain more precisely how this can be effectively realized
\emph{vis-a-vis} EM wave
propagation using metamaterials.

We first describe the wormhole as an abstract manifold $M$, see Fig. 
1; in the next section
we will show how to realize this concretely in $\R^3$, as a wormhole device
$N$. Start by making two holes in the Euclidian space
$\R^3=\{(x,y,z)|x,y,z\in\R\}$, say by removing the open ball $B_1=B({\it O},1)$
with center at the origin
${\it O}$ and of radius 1,
and also the open ball $B_2=B(P,1)$, where $P=(0,0,L)$ is a point on
the $z$-axis having the
distance $L>3$ to the origin. We denote  by $M_1$ the region so obtained,
$M_1=\R^3\setminus (B_1\cup B_2)$,
which is
the first component  we need to construct a
wormhole.
Note that $M_1$ is a 3-dimensional manifold with boundary, the boundary  of
$M_1$ being
$\p M_1=\p B_1\cup
\p B_2$,   the union of two 2-spheres. Thus, $\p M_1$
can be considered as a disjoint  union
$\sph\cup\sph$, where we will use $\S^2$ to denote various copies of the
two-dimensional unit sphere.

The second component needed is a $3-$dimensional cylinder, $M_2=\S^2
\times[0,1]$. This cylinder can be constructed by
taking the closed unit cube $[0,1]^3$ in $\R^3$ and,  for each value
of $0<s<1$,
gluing together, i.e., identifying, all of
the points on the boundary of the cube with $z=s$. Note that we do not identify
points at the top of the boundary, at  $z=1$, or at the bottom, at $z=0$.
We then glue together the boundary $\p B({\it O}, 1)\sim\S^2$ of the ball
$B({\it O}, 1)$
with the lower end (boundary component) $\S^2\times\{0\}$ of  $M_2$, and the
boundary
$\p B(P, 1)$ with the upper end, $\S^2\times\{1\}$.
In doing so we
identify the point $(0,0,1) \in  \p B({\it O}, 1)$ with the point $NP
\times \{0\}$
and the point $(0,0, L-1) \in  \p B(P, 1)$ with the point $NP \times \{1\}$,
where $NP$ is the north pole on $\S^2$.

The resulting manifold $M$  no longer lies in $\R^3$, but rather
    is the connected sum of the components $M_1$ and $M_2$ and has the topology
of $\R^3$ with a $3-$dimensional handle attached.
Note that adding this handle makes it possible to travel
        from one point in  $M_1$ to another point in  $M_1$, not only
        along curves lying in
       $M_1$ but also those in $M_2$.

To consider Maxwell's equations on
$M$,  let us start with
Maxwell's equations   on $\R^3$ at frequency $k\in\R$, given by
\ba
\nabla\times E = ik B,\quad \nabla\times H =-ikD,\quad
D(x)=\e(x) E(x),\quad  B(x)=\mu(x) H(x).
\ea
Here $E$ and $H$ are the electric and magnetic fields, $D$ and $B$
are the electric displacement field  and
the magnetic flux density,
$\e$ and $\mu$ are matrices corresponding to permittivity and permeability.
As the wormhole is topologically different from the
Euclidian space $\R^3$,
we  use  a formulation of Maxwell's equations on a manifold, and as in
\cite{GKLU}, do this in the setting of  a general Riemannian manifold,
$(M,g)$. For our
purposes, as in \cite{KLS,GKLU} it suffices to use
$\e, \mu$ which are conformal, i.e., proportional by scalar fields,
to the metric $g$.
In this case, Maxwell's equations  can be written, in the coordinate
invariant form,
as
\ba
& &d E = ikB,\ d H =-ikD,\quad
        D=\epsilon E,\  B=\mu H\quad \hbox{in }M,
\ea
where $E,H$ are 1-forms, $D,B$ are 2-forms, $d$ is the exterior derivative,
       and
$\epsilon$ and $\mu$ are scalar functions times the Hodge operator of
$(M,g)$, which maps 1-forms to the corresponding 2-forms \cite{Frankel}.
In  local coordinates these equations are written in the same
form as Maxwell's equations in Euclidian space with matrix valued
$\e$ and $\mu$.
Although not necessary,  for simplicity one can choose a metric on the wormhole
manifold $M$ which is
       Euclidian  on $M_1$, and  on $M_2$
       is the product of a given metric
$g_0$ on $\S^2$ and the standard metric of $[0,1]$.
More generally,  can also choose the metric on $M_2$ to  be a warped  product.
Even the simple choice
of the product of the standard metric of $\S^2$ and
the metric $\delta^2 ds^2$, where
$\delta$ is the ``length" of the wormhole,  gives rise to
interesting ray-tracing effects for rays passing through the wormhole
tunnel. For
$\delta <<1$,
    the image through one end of the wormhole  (of the region beyond the
other end) would resemble the
image in a a fisheye lens; for $\delta\gtrsim1$, multiple images and
greater distortion occur.
(See \cite[~Fig.2]{GKLU2}.)

The proof of the wormhole effect that we actually give is for yet another
variation, where the balls
that form the ends have their boundary spheres flattened;
this may be
useful for applications, since it allows for
there to be a vacuum (or
air) in a neighborhood of the axis of the wormhole, so that, e.g., instruments
may be passed through the wormhole.  We next  show how to construct, using
metamaterials, a device $N$  in
$\R^3$  that effectively realizes the  geometry and topology of $M$,
relative to solutions of Maxwell's equations at frequency $k$, and hence
functions as an electromagnetic wormhole.

\section{The wormhole device $N$ in $\R^3$}

We  now explain how to construct a  ``device" $N$ in $\R^3$,   i.e., a
specification of permittivity $\e$ and permeability $\mu$,  which  affects the
propagation of electromagnetic waves in the same way as the presence of the
handle $M_2$ in the wormhole manifold $M$. What this means is that we
prescribe
a configuration of  metamaterials which make the waves
behave as if there were an invisible tube attached to
$\R^3$, analogous to the handle $M_2$ in the wormhole manifold $M$.
In the other words, as far as external EM observations of the wormhole device
are concerned, it
appears  as if the topology of  space
has been  changed.

\begin{figure}[htbp]
\begin{center}
\psfrag{1}{$M$}
\psfrag{2}{$M_2$}
\psfrag{3}{$M_1$}
\includegraphics[width=8cm]{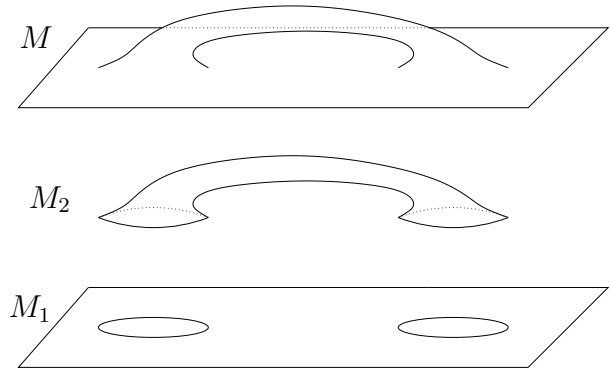} \label{pic 3}
\end{center}
\caption{Schematic figure: a
wormhole manifold is glued from two components,
the ``handle'' and space with two holes. Note that
in the actual construction, the components are three dimensional.}
\end{figure}

We use  cylindrical coordinates $(\theta,r,z)$ corresponding to  a point
$(r\cos\theta,r\sin\theta,z)$ in $\R^3$.
The wormhole device is built  around an obstacle $K\subset \R^3$.
To define $K$, let
$S$ be the  two-dimensional finite cylinder
       $\{\theta \in [0, 2\pi], r=2,\ 0\leq z\leq L\} \subset \R^3$. The
open region $K$
consists of
all points in $\R^3$ that have distance less than one to $S$ and has
the shape of  a long, thick-walled tube with smoothed corners.

Let us first introduce a deformation map $F$ from $M$ to
$N =\R^3\setminus K$ or,
more precisely, from $M \setminus \gamma$ to $N \setminus \Sigma$,
where $\gamma$ is a closed curve in $M$ to be described shortly and
$\Sigma = \p K$.
We will define $F$ separately on $M_1$ and $M_2$
denoting the corresponding parts by $F_1$ and $F_2$.

To describe $F_1$, let $\gamma_1$ be the line segment
on the $z-$axis connecting $\p B({\it O}, 1)$ and $\p B(P, 1)$ in $M_1$,
namely, $\gamma_1=\{r=0,\ z\in [1,L-1]\}$.
Let $
F_1(r,z)=(\theta,R(r,z),Z(r,z))$ be such that $(R(r,z),Z(r,z))$,
shown in  Fig.\ 2,

\begin{figure}[htbp]
\begin{center}
\psfrag{1}{$A$}
\psfrag{2}{$B$}
\psfrag{3}{$C$}
\psfrag{4}{$D$}
\psfrag{5}{$A'$}
\psfrag{6}{$B'$}
\psfrag{7}{$C'$}
\psfrag{8}{$D'$}
\psfrag{9}{}
\includegraphics[width=14cm]{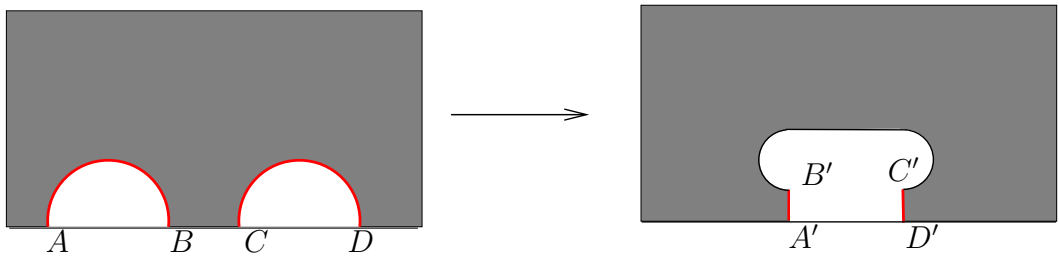} \label{pic 5}
\end{center}
\caption {The map $(R(r,z),Z(r,z))$ in  cylindrical coordinates  $(z,r)$.}
\end{figure}

transforms in the $(r,z)$ coordinates the
semicircles $AB$
and $CD$
in the left picture to the vertical line segments $A'B'=\{ r \in [0,
1], z=0\}$ and
$C'D'=\{ r \in [0, 1], z=L\}$ in the right picture and the cut $\gamma_1$
on the left picture to the curve $B'C'$ on the
right picture. This gives us a map
$
F_1:M_1\setminus \gamma_1 \to N_1 \setminus \Sigma,
$
where the closed region $N_1$ in $\R^3$ is obtained by  rotation of
the region exterior to the curve
$A'B'C'D'$ around the $z-$axis.
We can choose $F_1$ so that it is the identity map
in the domain $U=\R^3\setminus \{-2\leq z\leq L+2,\ 0\leq r\leq 4\}$.

To describe $F_2$,  consider the line segment,
$\gamma_2=\{NP\}\times [0,1]$   on $M_2$ . The sphere without the
north pole can
be "flattened"  and stretched to an open disc with radius one which,
together with stretching
$[0,1]$ to $[0, L]$,
       gives us a map
$F_2$ from $M_2\setminus \gamma_2$ to $N_2 \setminus \Sigma$. The region
$N_2$ is  the
$3-$dimensional cylinder, $N_2=\{\theta \in [0, 2\pi], r \in [0,
1], z \in [0,L]\}$.
When flattening $\S^2 \setminus NP$, we do it in such a way that
$F_1$ on $\p B({\it O}, 1)$ and $\p B(P, 1)$ coincides with $F_2$ on
$(\S^2 \setminus NP) \times \{0\}$ and $(\S^2 \setminus NP) \times \{1\}$,
respectively.

Thus, $F$ maps $M\setminus \gamma$, where
$\gamma=\gamma_1 \cup \gamma_2$ is a closed curve in $M$,
onto $N \setminus \Sigma$; in addition, $F$ is the identity
on the region $U$.

Now we are ready to define the electromagnetic  material parameter
tensors
on $N$. We define
the permittivity  to be
\ba
\tilde \e=F_*\e(y)=\left. \frac{(D F)(x)\cdotp \e(x)
\cdotp (D F(x))^t}{\det (D F)}\right|_{x=F^{-1}(y)},
\ea
where $DF$ is the derivative matrix of $F$, and similarly the
permeability to be
$\tilde \mu=F_*\mu$.
These deformation rules are based on the fact that permittivity
and permeability are  conductivity type tensors, see \cite{KLS}.

Maxwell's equations are invariant under smooth changes of coordinates.
This means that, by the chain rule, any solution to Maxwell's
equations in
$M\setminus \gamma$, endowed with material parameters $\e,\mu$
becomes, after  transformation by  $F$, a solution to
       Maxwell's equations in
$N \setminus \Sigma$ with material parameters $\tilde \e$ and $\tilde \mu$,
and {\it vice versa}. However, when considering the fields on the
entire spaces $M$
and
$N$, these observations are not enough, due to the singularities
of $\tilde \e$ and $\tilde \mu$ near $\Sigma$; the significance of
this for cloaking was observed and
analyzed in \cite{GKLU}. In the following,  we will show
that the physically
relevant class of solutions to  Maxwell's equations, namely the {\it
(locally) finite
energy}  solutions, remains the same, with respect to the transformation $F$,
in $(M; \e,\mu)$  and $(N; \tilde \e, \tilde \mu).$
One can analyze the rays  in
$M$ and $N$ endowed with the electromagnetic wave propagation
metrics $g = \sqrt{\e \mu}$ and $ \tilde g = \sqrt{ \tilde \e \tilde \mu}$,
respectively. Then the rays on $M$ are transformed by $F$ into the rays
in $N$. As almost all the rays on $M$ do not intersect with $\gamma$,
therefore, almost all the rays on $N$ do not approach $\Sigma$.
This was the basis for \cite{Le,PSS1} and was analyzed further in
\cite{PSS2}; see also \cite{MBW} for a similar analysis in the context of
elasticity.
Thus, heuristically one is  led to conclude that the
electromagnetic waves on
$(M; \e,\mu)$  do not feel the presence of $\gamma$, while those on
$(N; \tilde \e,
\tilde\mu)$ do not feel the presence  of $K$,
and these waves can be transformed into each other by the
map $F$.

Although the above considerations are  mathematically rigorous,
on the level both of the chain rule and  of high
frequency limits, i.e., ray tracing,  in the
exteriors $M\setminus \gamma$ and $N\setminus\Sigma$, they do not suffice to
fully describe the behavior of physically meaningful solution fields on $M$
and $N$.  However, by carefully examining the class of
the finite-energy waves
in $M$ and $N$ and analyzing their behavior near
$\gamma$ and $\Sigma$, respectively, we can give a complete analysis,
justifying the
conclusions above.
Let us briefly explain the main steps of the analysis
using methods developed for theory of invisibility (or cloaking)
     at frequency $k>0$
\cite{GKLU} and at frequency $k=0$ in \cite{GLU1,GLU2}.
The details will follow.

First, to guarantee that the fields in $N$ are finite energy
solutions and  do not blow up near  $\Sigma$, we have to
impose at $\Sigma$
the appropriate boundary condition, namely,
the  Soft-and-Hard (SH)  condition, see \cite{HLS,Ki2},
\ba
e_\theta\,\cdotp
E|_{\Sigma}=0,\quad
e_\theta\,\cdotp
H|_{\Sigma}=0,
\ea
where $e_\theta$ is the angular direction.
Secondly, the map $F$ can be considered as a smooth coordinate
transformation on $M\setminus\gamma$; thus, the
finite energy solutions on $M\setminus \gamma$
transform under  $F$ into  the
finite energy solutions on $N\setminus \Sigma$, and vice versa.
Thirdly,  the curve $\gamma$ in $M$ has  Hausdorff dimension equal to one.
This implies that the possible singularities of
the finite energy electromagnetic fields near $\gamma$
are removable \cite{KKM}, that is, the finite energy
fields in $M\setminus \gamma$ are exactly the restriction to
$M\setminus \gamma$ of the fields defined on all of $M$.

Combining these steps
we can see that measurements of the electromagnetic fields on
$(M;\e,\mu)$ and on $(\R^3\setminus K; \tilde \e, \tilde \mu)$
coincide in $U$.
In the other words, if we apply any current on $U$ and
measure the radiating electromagnetic fields it generates,
then the fields on $U$ in the wormhole manifold $(M; \e, \mu)$
coincide with the fields on $U$ in $(\R^3\setminus K; \tilde \e, \tilde \mu)$,
$3$-dimensional space equipped with the wormhole device construction.

Summarizing our construction, the wormhole device  consists of the
metamaterial coating of the obstacle
$K$. This coating should have the permittivity
       $\tilde \e$  and permeability $\tilde \mu$.
       In addition, we need to impose the SH boundary condition
on $\Sigma$, which may be realized
       by fabricating the obstacle  $K$ from a perfectly conducting 
material with
parallel corrugations on its surface \cite{HLS,Ki2}.

In the next section,  the
permittivity $\tilde \e$ and and permeability $\tilde \mu$
are described in a rather simple form.  (As mentioned earlier, in
order to allow for
a tube around the axis of the wormhole to be a vacuum or air, we deal with a
slightly different construction than was described above, starting
with flattened spheres).
It should be possible to physically implement an approximation to this
mathematical idealization of
the material parameters needed for the wormhole device,
using   concentric rings of split ring resonators as in
the experimental verification of cloaking
obtained in \cite{SMJCPSS}.

\section{Rigorous construction of the wormhole}

Here we present a rigorous model of a typical wormhole device and
justify the claims above concerning the behavior of the electromagnetic
fields in the wormhole
device in $\R^3$ in terms of  as the fields on the wormhole manifold $(M,g)$.

\subsection{The wormhole manifold $(M,g)$ and
the wormhole device $N$}

Here we prove the wormhole effect for a  variant of the wormhole
device described in the previous sections.
Instead of using a round sphere $\S^2$ as before, we present
a construction  that uses a deformed  sphere $\S_{\rm flat}^2$ that
is flat the near the south and north  poles, $SP$ and  $NP$.
This makes  it possible to have constant isotropic material parameters
near the $z$-axis located inside the wormhole. For possible applications, see
\cite{GKLU2}.

We use following notations.
Let $(\theta,r,z)\in  [0,2\pi]\times \overline  \R_+\times \R$ be
the cylindrical coordinates of $\R^3$, that is the map
\ba
X:(\theta,r,z)\to (r\cos\theta,r\sin\theta,z)
\ea
that maps $X:[0,2\pi]\times \overline \R_+\times  \R\to \R^3$.
In the following, we identify $[0,2\pi]$ and the unit circle $S^1$.

Let us start by removing from $\R^3$ two ``deformed" balls which have
flat portions near the south and north poles. More precisely, let
$M_1=\R^3\setminus (P_1\cup P_2)$, where in the cylindrical coordinates
\ba
P_1&=&\{X(\theta,r,z):\ -1\leq z\leq 1,\ 0\leq r\leq 1\}\\
& &\ \ \cup
\{X(\theta,r,z):\ (r-1)^2+z^2\leq 1\},\\
P_2&=&\{X(\theta,r,z):\-1\leq z-L\leq 1,\ 0\leq r\leq 1\}  \\
& &\ \ \cup
\{X(\theta,r,z):\ (r-1)^2+(z-L)^2\leq 1\}.
\ea
We say that the boundary $\p P_1$ of $P_1$ is a deformed sphere with
flat portions, and denote it by $\S^2_{\rm flat}$.
We say that the intersection  points of $\S^2_{\rm flat}$
with the $z$-axis are the north pole, $NP$, and the south pole,
$SP$.

Let $g_1$ be the metric on $M_1$ inherited from $\R^3$,
and let $\gamma_1$ be the path
\ba
\gamma_1=\{X(0,0,z):\ 1<z<L-1\}\subset M_1.
\ea
Set
\ba
A_1&=&M_1\setminus V_{1/4},\\
V_t&=&\{X(\theta,r,z):\ 0\leq r\leq t,\ 1< z< L-1\},\quad 0<t<1,
\ea
and consider a map $G_0:M_1\setminus \gamma_1\to A_1$; see Fig. 3. $G_0$
defined as
the identity map on $M_1\setminus V_{1/2}$ and,
in  cylindrical coordinates, as
\ba
G_0(X(\theta,r,z))= X(\theta,\frac 14+\frac r2,z),\quad
(\theta,r,z)\in V_{1/2}.
\ea
Clearly, $G_0$ is $C^{0,1}-$smooth.

Let $U(x)\in \R^{3\times 3}$, $x=X(\theta,r,z)$,
be the orthogonal matrix that maps
the standard unit vectors $e_1,e_2,e_3$ of $\R^3$ to the Euclidian unit vectors
corresponding to the $\theta$, $r$, and $z$ directions, that is,
\ba
U(x)e_1=(-\sin \theta,\cos \theta,0),\quad
U(x)e_2=(\cos \theta,\sin \theta,0),\quad
U(x)e_3=(0,0,1).
\ea
Then the differential of $G_0$ in the Euclidian
coordinates at the point $x\in V_{1/2}$
      is the matrix
\beq\label{eq: diff G_0}
DG_0(x)
U(y)\left(\begin{array}{ccc} \frac 1r(\frac 14+\frac r2) &0 &0 \\
0 &\frac 12 &0 \\
0 &0 &1 \end{array}\right)U(x)^{-1},\quad x=X(\theta,r,z),\ y=G_0(x).
       \eeq
Later we impose on  part of the boundary, $\Sigma_0=\p A_1\cap
\{1<z<L-1\}$,
the soft-and-hard boundary condition (marked red
in the figures).

Next, let
$(\theta,z,\tau)=(\theta(x),z(x),\tau(x))$ be the Euclidian boundary
normal coordinates
associated to $\Sigma_0$, that is, $\tau(x)=\dist_{\R^3}(x,\Sigma_0)$
and $(\theta(x),z(x))$ are the $\theta$ and $z$-coordinates
of the closest point of  $\Sigma_0$ to $x$.

Denote by $(G_0)_*g_1$
the push forward of the metric $g_1$ in $G_0$, that is, the metric obtained
from $g_1$ using the change of coordinates $G_0$, see \cite{Frankel}.
The metric $(G_0)_*g_1$ coincides
with $g_1$ in $A_1\setminus V_{1/2}$,
and in the Euclidian boundary normal coordinates
of $\Sigma_0$,  on
$A_1\cap V_{1/2}$, the metric  $(G_0)_*g_1$,
     has  the length element
\ba
ds^2=4\tau^2\,d\theta^2+
dz^2+4d\tau^2.
\ea

\begin{figure}[htbp]
\begin{center}
\psfrag{1}{$G_0$}
\includegraphics[width=12cm]{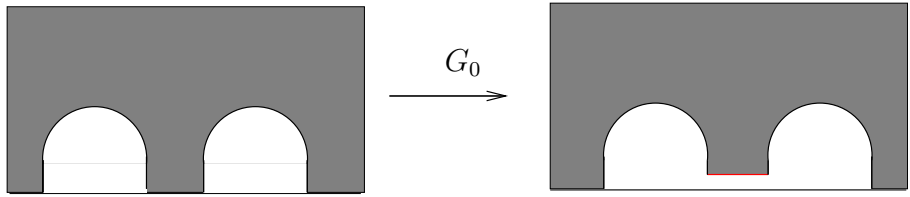} \label{pic F1}
\end{center}
\caption{A schematic figure on the map $G_0$,
considered in the $(r,z)$ coordinates.
Later, we impose the SH boundary condition on the portion of the boundary
coloured red.}
\end{figure}

Next, let
\ba
q_3&=&\hbox{conv}\bigg(
\{(r,z) :\ (r-2)^2+(z-(-2))^2\leq 1\}
\\
& & \quad \quad \quad \quad \cup
\{(r,z) :\ (r-2)^2+(z-(L+2))^2\leq 1\}\bigg),\\
q_4&=&\{(r,z) :\ 0\leq r \leq 1,\ -1\leq z\leq L+1\},
\ea
where conv$(q)$ denotes the convex hull of the set $q$.

Let
\ba
N_1&=&\R^3\setminus (P_3\cup P_4),\\
P_3&=&\{X(\theta,r,z):\ (r,z)\in q_3\},\\
P_4&=&\{X(\theta,r,z):\ (r,z)\in q_4\},\\
\Sigma_1&=&\p N_1\setminus \p P_4.
\ea
We can find a   Lipschitz smooth map $G_1:A_1\to N_1$,
see Fig. 4,  of the form
\ba
G_1(X(\theta,r,z))=X(\theta,R(r,z),Z(r,z))
\ea
such that it maps
$\Sigma_0$ to $\Sigma_1$, and in $A_1$ near $\Sigma_0$
it is given by
\beq\label{eq: normal 1}
G_1(x+t\nu_0)=G_1(x)+t\nu_1.
\eeq
Here, $x\in \Sigma_0$, $\nu_0$ is the Euclidian unit
normal vector of $\Sigma_0$, $\nu_1$ is the Euclidian unit
normal vector of $\Sigma_1$, and $0<t<\frac 14$. Moreover,
we can find a $G_1$ so that it is the identity map near the $z$-axis,
that is,
\beq\label{eq: normal 2}
G_1(x)=x,\quad x\in A_1\cap \{0\leq r< \frac 14\}
\eeq
and such that  $G_1$ is also the identity map in the set of points with
the Euclidian distance 4 or more
from $P_1\cup P_2$.
Note that we can find such a $G_1$ such that both $G_1$ and its inverse
$G_1^{-1}$ are   Lipschitz smooth  up to the boundary. Thus
the differential $DG_1$ of $G_1$ at $x\in A_1$ in
Euclidian coordinates is
\ba
DG_1(x)=
U(y)\left(\begin{array}{cc} a_{11}(r,z) &0 \\
0 & A(r,z)\end{array}\right)U(x)^{-1},\quad x=X(\theta,r,z),\ y=G_1(x),
\ea
where $c_0\leq a_{11}(r,z)\leq c_1$
and $A(r,z)$ is a symmetric $(2\times 2)$-matrix
satisfying
\ba
c_0 I\leq A(r,z)\leq c_1I
\ea with some
$c_0,c_1>0$.

{The map $F_1(x)=G_1(G_0(x))$ then maps
$F_1:M_1\setminus \gamma_1\to N_1$. } Let $\tilde g_1=(F_1)_*g_1$
be metric on $N_1$.
   From the above considerations, we see
      that
the differential $DF_1$ of $F_1$ at $x\in M_1\setminus \gamma_1$
near $\Sigma_0$, in
Euclidian coordinates, is given by
\beq\label {eq: diff F_1}
& &DF_1(x)=
U(y)\left(\begin{array}{cc} b_{11}(\theta,r,z) &0 \\
0 & B(r,z)\end{array}\right)U(x)^{-1},
\\
\nonumber
& & b_{11}(\theta,r,z)=\frac
{c_{11}(r,z)}{\dist_{\R^3}(X(\theta,r,z),\Sigma_0)}\quad 
x=X(\theta,r,z),\ y=F_1(x)
  \eeq
where   $c_0\leq c_{11}(r,z)\leq c_1$,
and $B(r,z)$ is a symmetric $(2\times 2)$-matrix
satisfying
\ba
c_0 I\leq B(r,z)\leq c_1I,
\ea
for some
$c_0,c_1>0$.

Note that
$\p P_4\cap \{r<1\}$ consists of two two-dimensional discs,
$B_2(0,1)\times \{-1\}$ and $B_2(0,1)\times \{L+1\}$.
Below, we will
use the map
\ba
f_2=F_1|_{\p P_1\setminus NP}:{\p P_1}\setminus NP\to
B_2(0,1)\times \{-1\}\subset \p N_1.
\ea
The map $f_2$ can be considered as the deformation
that ``flattens'' $\S^2_{\rm flat}\setminus NP$ to a two
dimensional unit disc.

\begin{figure}[htbp]
\begin{center}
\psfrag{1}{$G_1$}
\includegraphics[width=12cm]{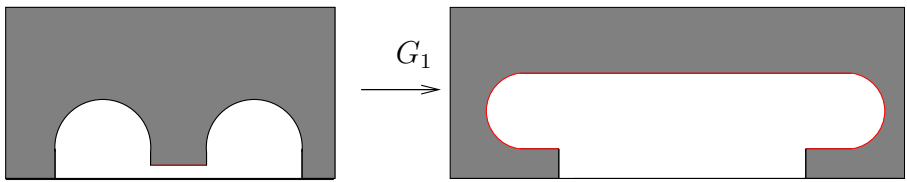} \label{pic F3}
\end{center}
\caption{Map $G_1$ in $(r,z)$-coordinates.}
\end{figure}

To describe $f_2$, consider $\S^2_{\rm flat}$ as a surface in  Euclidian space
and define on it the $\theta$ coordinate corresponding
to the $\theta$ coordinate of $\R^3\setminus \{z=0\}$.
Let then $s(y)$ be
the intrinsic distance of $y\in\S^2_{\rm flat}$ to the south pole $SP$.
Then  $(\theta,s)$ define coordinates in $\S^2_{\rm flat}\setminus
\{SP,NP\}$.
We denote by $y(\theta,s)\in \S^2_{\rm flat}\setminus
\{SP,NP\}$ the point  corresponding to  the coordinates  $(\theta,s)$.

By the above construction, the map $f_2$ has  the form,
     {with respect to the coordinates used
above, }
\beq\label{eq: f_2 structure}
& &f_2(y(\theta,s))=X(\theta,R(s),-1)\in B_2(0,1)\times \{-1\},\ \ \
\hbox{where }\\ \nonumber
& &R(s)=s,\quad \hbox{for }0<s<\frac 14,
\\ \nonumber
& &R(s)=1-\frac12[(\pi+4)-s],\quad \hbox{for }(\pi+4)-\frac 14 <s<(\pi+4),
\eeq
cf.\ formulae (\ref{eq: normal 1}) and (\ref{eq: normal 2}).
In the following  we identify $B_2(0,1)\times \{-1\}$ with the disc $B_2(0,1)$.

Let $h_1$ be the metric on $\p P_1\setminus NP$ inherited
from $(M_1,g_1)$. Let $h_2=(f_2)_*h_1$ be the metric on
$B_2(0,1)$.
We observe  that
the metric $h_2$ makes the disc $B_2(0,1)$
isometric to $\S^2_{\rm flat}\setminus NP$,
endowed with the metric inherited from $\R^3$. Thus, let
\ba
M_2=\S^2_{\rm flat}\times [-1,L+1].
\ea
      On $M_2$, let the metric $g_2$ be
the product of the metric of $\S^2_{\rm flat}$ inherited
from $\R^3$ and the metric $\a_2(z) dz^2,\, \a_2 >0$ on $[-1, L+1]$.
Let $\gamma_2=\{NP\}\times [-1,L+1]$ be a path on $M_2$.

Define $N_2=P_4=\{X(\theta,r,z):\ 0\leq r<1, -1\leq z\leq L+1\} \subset \R^3$,
$\Sigma_2=\p N_2\cap \{r=1\}$,
    and
let $F_2:M_2\setminus \gamma_2\to N_2$ be the map
of the form
\beq\label{eq: F_2 structure}
F_2(y,z)=(f_2(y),z)\in \R^3,\quad (y,z)\in(\S^2_{\rm flat}\setminus
NP)\times [-1,L+1].
\eeq
Let $\tilde g_2=(F_2)_*g_2$ be  the resulting metric on $N_2$.

\begin{figure}[htbp]
\begin{center}
\includegraphics[width=6cm]{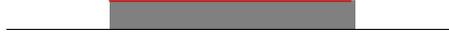} \label{pic F2}
\end{center}
\caption{The set $N_2$ in the $(r,z)$ coordinates.
Later, we  impose the SH boundary condition on the portion of the boundary
colored red.}
\end{figure}

Denote by $\overline M_1=M_1\cup \p M_1$  the closure
of $M_1$ and let $(M,g)=(\overline M_1,g_1)\# (M_2,g_2)$
be the connected sum of
$\overline M_1$ and $M_2$, that is, we glue
the boundaries $\p M_1$ and $\p M_2$.
The set $N=N_1\cup N_2\subset \R^3$ is open, and
its boundary $\p N$ is $\Sigma=\Sigma_1\cup \Sigma_2$.

Let $F$ be the map
$F:M\setminus \gamma \to N$ defined by
the maps $F_1:M_1\setminus \gamma_1\to N_1$ and  $F_2:M_2\setminus
\gamma_2\to N_2$, and finally,
let $\gamma=\gamma_1\cup \gamma_2$ and
$\tilde g=F_*g$.

\begin{figure}[htbp]
\begin{center}
\includegraphics[width=6cm]{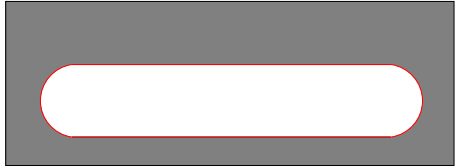} \label{pic F3B}
\end{center}
\caption{The set $N=N_1\cup N_2\subset \R^3$ having the complement
$K$,
presented in the $(r,z)$ coordinates.
Later, the SH boundary condition is imposed on $\p K$.}
\end{figure}

Let $K=\R^3\setminus N$. On the surface
$\Sigma=\p K$ we can use local coordinates $(\tilde t,\tilde \theta)$, where
$\tilde \theta$ is the $\theta$-coordinate of the ambient space $\R^3$
and $\tilde t$ is either the $r$ or $z$ -coordinate of
the  ambient space $\R^3$ {\newtekst restricted to $\Sigma$}. Denote also
\ba
\tilde \tau=\tilde \tau(x)=\dist_{\R^3}(x,\p K).
\ea
    Then by
formula (\ref{eq: normal 1})
we see that in $N_1$, in the Euclidian boundary normal
coordinates $(\tilde \theta,\tilde t,\tilde \tau)$ associated to the surface
$\Sigma_1$,
the metric $\tilde g$ has the length element
\ba
ds^2=4d\tilde \tau^2+\a_1(\tilde t)\, d\tilde t^2+4\tilde \tau^2\,d
\tilde \theta^2,\quad 0<\tilde \tau<\frac 14,\quad
c_0^{-1} \leq \a_1(\tilde t)\leq c_0,\, c_0\geq 1.
\ea
The construction of $F_2$ yields that in $N_2$ ,
      in the Euclidian boundary normal
coordinates $(\tilde \theta,\tilde t,\tilde \tau)$
with $\tilde t=z$, associated to the surface
$\Sigma_2=\p K\cap \p N_2$, the metric $\tilde g$ has the length element,
{\newtekst near $\Sigma_2$},
\ba
ds^2=4d\tilde \tau^2+\alpha_2(\tilde t)d\tilde t^2+4\tilde
\tau^2\,d\tilde
\theta^2,\quad 0<
\tilde \tau<\frac 14.
\ea
Here,  near $\p N_1\cap \p N_2$,  we use $\tilde t=z$ on  $\Sigma_1$.
Choosing the map $G_1$ in the construction of the map $F_1$
appropriately, we have $\a_2(-1)=\a_1(-1),\,
\a_2(L+1)=\a_1(L+1)$, and  the resulting map is Lipschitz.

On $M_1$, $N_1$, and $N_2$ that are subsets of $\R^3$ we have
the well defined cylindrical coordinates $(\theta,r,z)$.
Similarly, $M_2=\S^2_{\rm flat}\times [-1,L+1]$ we define
the coordinates $(\theta,s,z)$, where $(\theta,s)$ are
the above defined
coordinates on $\S^2_{\rm flat}\setminus
\{SP,NP\}$.

We can also consider on $N\subset \R^3$ also
the Euclidian metric, denoted by $g^e.$ In Euclidean
coordinates, $(g^e)_{ij}=\delta_{jk}$. Consider next
the above defined Euclidian boundary normal coordinates
$(\tilde \theta,\tilde t,\tilde \tau)$ associated to $\p K$.
They are well defined in a neighborhood of $\p K$.
We define the vector fields
\ba
\tilde \xi=\p_{\tilde \tau},\quad
\tilde \eta=\p_{\tilde \theta},\quad
\tilde \zeta=\p_{\tilde t}
\ea
on $N$ near $\p K$. These vector fields are orthogonal
with respect to the metric $\tilde g$ and to the metric $g^e$.

On $M$ near $\gamma$, we use coordinates
$(\theta,t,\tau)$.
On $M_1$, near $\gamma_1$ they
in the terms of the cylindrical coordinates are
$(\theta,t,\tau)=(\theta,z,r)$. On
$M_2$, they are the coordinates
$(\theta,t,\tau)=(\theta,z,s)$, where $s$
is the intrinsic distance to the north pole $NP$.
We define also the vector fields
\ba
\xi=\p_{\tau},\quad
\eta=\p_{\theta},\quad
\zeta=\p_{t}
\ea
on $M\setminus \gamma$ near $\gamma$. These vector fields are orthogonal
with respect to the metric \nolinebreak$g$.

In the sequel, we consider
the differential of $F$ as the linear map \linebreak$DF:(T_xM,g)\to
(T_yN,g^e)$, $y=F(x)$, $x\in M\setminus \gamma$.

Using formula (\ref{eq: diff F_1}) in $M_1$ and formulas (\ref{eq:
f_2 structure}),
(\ref{eq: F_2 structure}) in $M_2$,  we see that
$DF^{-1}(x)$ at $x\in N$ near $\p N$
is a bounded linear map that
satisfies
\beq\label{eq: eta and zeta estimates} \\
\nonumber
|(\eta,DF^{-1}(x)\tilde \eta)_g|\leq C\,\tilde \tau(x),\quad &
(\zeta,DF^{-1}(x)\tilde \eta)_g=0,&\quad
(\xi,DF^{-1}(x)\tilde \eta)_g=0,\\
\nonumber
(\eta, DF^{-1}(x)\tilde \zeta)_g=0,\quad&
|(\zeta, DF^{-1}(x)\tilde \zeta)_g|\leq C,&\quad
|(\xi, DF^{-1}(x)\tilde \zeta)_g|\leq C,\\
\nonumber
(\eta, DF^{-1}(x)\tilde \xi)_g=0,\quad &
|(\zeta, DF^{-1}(x)\tilde \xi)_g|\leq C,&\quad
|(\xi, DF^{-1}(x)\tilde \xi)_g|\leq C,
\eeq
where $C>0$  and $(\cdotp\, ,\cdotp)_g$ is the inner product
defined by the metric $g$. Moreover, we obtain similar estimates for $DF$
in terms of the Euclidian metric $g^e$,
\beq\label{eq: eta and zeta estimates 2} \\
\nonumber
|(\tilde\eta,DF(y) \eta)_{g^e}|\leq C\,\tau(y)^{-1},\quad &
(\tilde \zeta,DF(y)\eta)_{g^e}=0,&\quad
(\tilde \xi,DF(y)\eta)_{g^e}=0,\\
\nonumber
(\tilde \eta, DF(y)\zeta)_{g^e}=0,\quad&
|(\tilde \zeta, DF(y)\zeta)_{g^e}|\leq C,&\quad
|(\tilde \xi, DF(y)\zeta)_{g^e}|\leq C,\\
\nonumber
(\tilde \eta, DF(y) \xi)_{g^e}=0,\quad &
|(\tilde \zeta, DF(y) \xi)_{g^e}|\leq C,&\quad
|(\tilde \xi, DF(y) \xi)_{g^e}|\leq C
\eeq
for $y\in M\setminus \gamma$ near $\gamma$ with $C>0$.

Next, consider  $DF(y)$ at $y\in M\setminus \gamma$. Recall that the
singular values $s_j(y)$, $j=1,2,3$
      of $DF(y)$
are  the square roots of the
eigenvalues of $(DF(y))^tDF(y)$,
where  $(DF)^t$ is the transpose of $DF$.
By  (\ref{eq: eta and zeta estimates}),
the singular values $s_j=s_j(y)$, $j=1,2,3$,
      of $DF(y)$, numbered in  increasing order,
satisfy
\ba
& &c_1\leq s_1(y)\leq  c_2,\\
& &c_1\leq s_2(y)\leq  c_2,\\
& &\frac {c_1}{\tau(y)}\leq s_3(y)\leq  \frac {c_2}{\tau(y)},
\ea
where $c_1,c_2>0$.

The determinant of the matrix $DF(y)$ can be computed
in terms of its singular values by
$\det(DF)=s_1s_2s_3$. Later, we need
the norm of the matrix $\det(DF(y))^{-1}\, DF(y)$. It
satisfies by formula (\ref{eq: eta and zeta estimates 2})
\beq\label{eq: lemma identity}
\|\det(DF(y))^{-1}\, DF(y)\|=\|(\prod_{k=1}^3 s_k^{-1})
\diag(s_1,s_2,s_3)
\|=\max_{1\leq j\leq 3}\prod_{k\not= j} s_k^{-1}\leq c_1^{-2}.
\eeq

\subsection{Maxwell's equations on the wormhole with
SH coating }\label{four-two}

Let $dV_0(x)$ denote the Euclidian volume element on $\overline N\subset \R^3$.
Recall that $N\subset \R^3$ is open set with boundary $\p N=\Sigma$.
Let $dV_g$ be the Riemannian volume on $(M,g)$.
We consider below the map $F:M\setminus \gamma\to N$
as a coordinate deformation. The map $F$ induces for any differential
form $\tilde E$ on $N$ a form $E=F^*\tilde E$ in $M\setminus \gamma$
called the pull back of $\tilde E$ in $F$, see \cite{Frankel}.

Next, we consider Maxwell equations with degenerate material parameters
$\tilde \e$ and $\tilde \mu$ on $N$ with SH boundary conditions on
$\Sigma$. On $M$ and $N$ we define the
permittivity
and permeability by
setting
\beq \label{e and mu definitions}
&
&\e^{jk}=\mu^{jk}=\det(g)^{1/2}g^{jk},\quad \hbox{on }M,\\ \nonumber
& &\tilde \e^{jk}=\tilde \mu^{jk}=\det(\tilde g)^{1/2}\tilde g^{jk},
\quad \hbox{on }N.
\eeq
Here, and below, the matrix $[g_{jk}(x)]$ is the representation of the
metric $g$
in local coordinates, $[g^{jk}(x)]$ is the inverse
of the matrix $[g_{jk}(x)]$, and $\det(g)$ is the
determinant of $[g_{jk}(x)]$.
We note that the metric $\tilde g$ is degenerate near $\Sigma$,
and thus $\tilde \e$ and $\tilde \mu$, represented
as matrices in the Euclidian coordinates, have elements
that tend to infinity at $\Sigma$, that is,
the matrices $\tilde \e$ and $\tilde \mu$ have a singularity
near $\Sigma$.
\medskip

\noindent{\bf Remark 1.} Modifying the above construction
by replacing $M_2$ with $M_2=\S^2_{\rm flat}\times [l_1,l_2]$
for appropriate $l_1,l_2\in \R$ and choosing $F_1$ in an
appropriate way,
we can use local coordinates $(\tilde \theta,\tilde t)$ on $\Sigma$
  such that
the Euclidian distance along $\Sigma$ of points
$(\tilde \theta,\tilde t_1)$ and $(\tilde \theta,\tilde t_2)$
  is proportional to  $|\tilde t_1-\tilde t_2|$, and the
metric $\tilde g$  in
the Euclidian boundary normal coordinates $(\tilde \theta,
\tilde t,\tilde \tau)$
associated to $\p K$ has the form
\ba
ds^2=4d\tilde \tau^2+ d\tilde t^2+4 \tilde \tau^2\,d\tilde \theta^2,
\quad 0<\tilde \tau<\frac 14.
\ea
The metric corresponding to the
metamaterials used in the physical experiment in \cite{SMJCPSS} has the
same form in
Euclidian
boundary normal coordinates associated to an infinitely long cylinder
$B_2(0,1)\times \R$.
Thus it seems likely that metamaterials similar to
those used in the  experimental verification of cloaking
could be used to create physical wormhole devices working at  microwave
frequencies.

\subsection{Finite energy solutions  of Maxwell's equations
   \\ and the equivalence theorem}

In the following, we consider 1-forms $\tilde E=\sum_j\tilde E_jd\tilde x^j$
and  $\tilde H=\sum_j \tilde H_jd\tilde x^j$ in the Euclidian
coordinates  $(\tilde x^1,\tilde x^2,\tilde x^3)$ of $N\subset \R^3$.
In the sequel, we use Einstein's summation convention and omit the
sum signs. We use the Euclidian coordinates as we want to consider $N$ with
the differential structure inherited from the Euclidian space.
We say that $\tilde E_j$ and $\tilde H_j$ are the (Euclidian) coefficients
of the forms $\tilde E$ and $\tilde H$, correspondingly.
We say that these coefficients are in $L^p_{loc}(\overline N,dV_0)$,
$1\leq p<\infty$,
if
\ba
\int_W |E_j(x)|^p\,dV_0(x)<\infty,\quad\hbox{for all
bounded measurable sets}\ \
W\subset \overline N.
\ea

\begin{definition}\label{SHS-def}
We say that the 1-forms $ \tilde E$ and $ \tilde H$ are \emph{finite energy}
solutions of  Maxwell's equations in $N$
with the soft-and-hard (SH) boundary conditions
on $\Sigma$ and the frequency $k\not =0$,
\ba
& &\nabla\times \tilde E = ik
\tilde \mu(x)  \tilde H,\quad \nabla\times
                \tilde H =-ik  \tilde
\e(x) \tilde E+\tilde J\quad
\hbox{ on }N,\\
& &\tilde \eta\,\cdotp \tilde
E|_{\Sigma }=0,\quad
\tilde \eta\,\cdotp \tilde
H|_{\Sigma}=0,
\ea
if 1-forms $ \tilde E$ and $ \tilde H$
and 2-forms $\tilde D=\tilde\e\tilde E$ and $\tilde B=\tilde\mu \tilde H$
     in $N$
have coefficients in $L^1_{\rm loc}(\overline N,dV_0)$ and
satisfy
\ba
\|\tilde E\|_{L^2(W,|\tilde g|^{1/2}dV_0)}^2=
\int_{W} \tilde \e^{jk}\, \tilde E_j\, \overline{\tilde E_k}
\,dV_0(x)<\infty,\\
\|\tilde H\|_{L^2(W,|\tilde g|^{1/2}dV_0))}^2=\int_{W}
\tilde \mu^{jk}\, \tilde H_j\, \overline{\tilde H_k} \,dV_0(x)<\infty
\ea
for all bounded  measurable sets $W\subset N$,
and
finally,
\ba
& &\int_{N} ((\nabla\times \tilde
h)\,\cdotp \tilde E-
                 ik \tilde h \,\cdotp \tilde \mu(x)\tilde H)
\,dV_0(x)=0,\\
& &\int_N (
(\nabla\times \tilde e)\,\cdotp \tilde H+
\tilde e \,\cdotp ( ik\tilde \e(x)\tilde E-\tilde J))
\,dV_0(x)=0,
\ea
for all 1-forms $\tilde e$ and $\tilde h$ with coefficients in
$C^\infty_0(\overline N)$
that satisfy
\beq\label{special electric and magnetic condition 2 A1}
&
&\tilde \eta\,\cdotp \tilde e|_{\Sigma}=0,\quad \tilde \eta\,\cdotp \tilde
h|_{\Sigma}=0,
\eeq
where $\tilde \eta=\p_\theta$ is the angular vector field that is tangential to
$\Sigma$.
\end{definition}

Below, we use for 1-forms $E=E_jdx^j$ and $H=H_jdx^j$, given in
     local coordinates
$(x^1,x^2,x^3)$ on $M$, the notations
\ba
\nabla\times E=dH,\quad \nabla\cdotp(\e E)=d*E,\quad
      \nabla\cdotp(\mu H)=d*H,
\ea
where $d$ is the exterior derivative and $*$ is the Hodge operator
on $(M,g)$,
cf.\ formula (\ref{e and mu definitions}).

We have the following ``equivalent behavior of electromagnetic
fields on $N$ and $M$'' result,  analogous to the results
of \cite{GKLU} for cloaking.

\begin{theorem}\label{Equivalence thm}
Let $E$ and $H$ be 1-forms
on $M\setminus \gamma$ and $\tilde E$ and $\tilde H$ be
1-forms
with coefficients in $L^1_{\rm loc}(\overline N,dV_0)$
such that $E=F^*\tilde E$, $H=F^*\tilde H$.
Let $\tilde J$ and $J=F^* \tilde J$ be  2-forms
with smooth coefficients
in $N$ and $M\setminus \gamma$
that are supported away
from $\Sigma$ and  $\gamma$.

Then the following are equivalent:
\begin{enumerate}

\item
On $N$, the
1-forms $ \tilde E$ and $ \tilde H$ satisfy Maxwell's equations
with SH boundary conditions
in the sense of Definition \ref{SHS-def}.

\item
On $M$, the forms $E$ and $H$ can be extended on $M$ so that
they are classical solutions $E$ and $H$
of Maxwell's equations,
\ba
& &\nabla\times     E = ik    \mu   H,\quad
\hbox{in }M,\\
\nonumber
& &\nabla\times H =-ik  \e   E + J,\quad
\hbox{in }M.
\ea
\end{enumerate}
\end{theorem}

\noindent
{\bf Proof.}
Assume first that $E$ and $H$ satisfy Maxwell's equations
on $M$ with source $J$ supported away from $\gamma$.
Then $E$ and $H$ are $C^\infty$ smooth near $\gamma$.

Using $F^{-1}:N\to M\setminus \gamma$ we define the
      1-forms $\tilde E,\tilde H$  and 2-form $\tilde J$
on $N$ by
$\tilde E=(F^{-1})^* E$, $\tilde H=(F^{-1})^*
H$, and $\tilde J=(F^{-1})^* J.$
These fields satisfy  Maxwell's equations
in $N$,
\beq\label{eq: physical Max2 B}
\nabla\times
\tilde E = ik \tilde \mu(x)  \tilde H,\quad \nabla\times
                \tilde H
=-ik  \tilde \e(x) \tilde E+\tilde J\quad
\hbox{ in }N.
\eeq
Now, writing $E=E_j(x)dx^j$ on $M$ near $\gamma$, we see
using
the transformation rule for differential
1-forms that
the form
$\tilde E=(F^{-1})^*E$ is in local
coordinates
\beq\label{eq: (B)}
\tilde E=\tilde
E_j(\tilde x)d\tilde x^j,\quad
\tilde E_j(\tilde x)
=(DF^{-1})_j^k(\tilde
x)\,E_k(F^{-1}(\tilde
x))
,\quad
\tilde x\in N.
\eeq
Using the smoothness of $E$ and $H$ near $\gamma$ on $M$ and
formulae (\ref{eq: eta and zeta estimates}), we see
that
$\tilde E$, $\tilde H$ are forms on $N$ with $L^1_{\rm loc}(\overline
N,dV_0)$ coefficients.
Moreover,
\ba
\tilde \e(x)\tilde E(x)&=&
\det(DF(y))^{-1}
DF(y)\e(y)DF(y)^t (DF(y)^t)^{-1}E(y)\\
&=&\det(DF(y))^{-1}
DF(y)\e(y)E(y)
\ea
where $x\in N$, $y=F^{-1}(x)\in M\setminus \gamma$.
     Formula (\ref{eq: lemma identity}) shows that
$\tilde D=\tilde \e\tilde E$, and $\tilde B=\tilde \mu\tilde H$
are 2-forms on $N$ with $L^1_{\rm loc}(\overline N,dV_0)$ coefficients.

Let $\Sigma(t)\subset \overline N$ be the
$t$-neighbourhood
of
$\Sigma$ in the $\tilde g$-metric.
Note that for small $t>0$ the set
$\Sigma(t)$ is the Euclidian $(t/2)$-neighborhood of $\p K$.
Denote by
$\nu$ be the unit exterior Euclidian
normal vector of
$\p\Sigma(t)$ and the Euclidian inner product by
$(\tilde \eta,\tilde E)_{g^e}=\tilde \eta\,\cdotp\tilde E.$

Formulas (\ref{eq: eta and zeta estimates})  and (\ref{eq: (B)}) imply
      that the angular components satisfy
       \ba
|\tilde \eta\,\cdotp \tilde E|\leq Ct,\quad x\in
\p\Sigma(t),
\ea
and
\ba
|\tilde \zeta\,\cdotp \tilde E|\leq C,\quad x\in
\p\Sigma(t)
\ea
with some $C>0$.
Thus denoting by $dS$ the Euclidian surface area on $\p \Sigma(t)$,
Stokes' formula, formula (\ref{eq: physical Max2 B}),
     and the identity $\nu\times \tilde \xi=\pm\tilde \eta$
yield
\ba
&
&\int_{N}
((\nabla\times \tilde h)\,\cdotp \tilde E-
                 ik \tilde h
\,\cdotp \tilde \mu \tilde H) \,dV_0(x)\\
& &=\lim_{t\to 0}
\int_{{N}\setminus \Sigma(t)}
((\nabla\times \tilde h)\,\cdotp \tilde E-
ik \tilde h \,\cdotp \tilde \mu \tilde H) \,dV_0(x)
\nonumber
\\
& &=-\lim_{t\to 0} \int_{\p \Sigma(t)}
(\nu \times \tilde E)\,\cdotp \tilde h\,  dS(x) \nonumber
\\
& &=-\lim_{t\to 0} \int_{\p \Sigma(t)}
\nu \times  ((\tilde \eta\,\cdotp \tilde E)\tilde \eta +(
\tilde\zeta\,\cdotp \tilde E)\tilde \zeta )
\,\cdotp \tilde h\,  dS(x) \nonumber
\\ \nonumber
& &=0
\ea
for a test
function $\tilde h$ satisfying formula
(\ref{special electric and magnetic condition 2 A1}).

Similar analysis for $\tilde H$ shows
that
1-forms $ \tilde E$ and $ \tilde H$ satisfy Maxwell's
equations
with SH boundary conditions
in the sense of Definition
\ref{SHS-def}.

Next, assume that $\tilde E$ and $\tilde H$ form a
finite energy solution of
Maxwell's equations
on $(N,g)$
with a source $\tilde J$
supported
away from $\Sigma$, implying in
particular that
\ba
\tilde
\e^{jk}\tilde E_j\overline {\tilde
E_k} \in L^1(W,\,dV_0),\quad
\tilde \mu^{jk}\tilde H_j\overline
{\tilde H_k} \in
L^1(W,\,dV_0)
\ea
where $W=F(U\setminus \gamma) \subset N$ and
$U\subset M$ is a relatively compact open neighbourhood of $\gamma$,
$\supp(\tilde J)\cap W=\emptyset$.
Define  $E=F^*\tilde E$,
$H=F^*\tilde H$, and
$J=F^*\tilde J$ on $M\setminus \gamma$. Therefore we conclude that
\ba
\nabla\times E =
ik \mu(x) H,\quad \nabla\times H =-ik \e(x)
E+J,
\quad \hbox{in }
M\setminus \gamma
\ea
and
\ba
\e^{jk}  E_j \overline {E_k} \in
L^1(U\setminus \gamma,\,dV_g),\quad
\mu^{jk} H_j  \overline {H_k} \in L^1(U\setminus \gamma,\,dV_g).
\ea
As representations of $\e$ and $\mu$, in local coordinates of
$M$, are matrices that are bounded from above and below, these
imply that
\ba
& &
\nabla\times E \in L^2(U\setminus \gamma,\,dV_g),\quad
\nabla\times H \in L^2(U\setminus \gamma,\,dV_g),\\
& &\nabla\cdotp (\e E)=0,\quad
\nabla\cdotp (\mu H)=0,\quad \hbox{in } U\setminus \gamma.
\ea
Let
$E^e,H^e\in L^2(U,\,dV_g)$ be measurable extensions
of $E$ and $H$
to $\gamma$.
Then
\ba
                & &\nabla\times E^e
-ik \mu(x) H^e=0,\quad \hbox{in }
U\setminus \gamma,\\
                &
&\nabla\times E^e -ik \mu(x) H^e\in H^{-1}(U,\,dV_g),\\
       & &\nabla\times H^e +ik \e(x) E^e=0,\quad \hbox{in }
U\setminus
\gamma,\\
                & &\nabla\times H^e +ik \e(x) E^e\in
H^{-1}(U,\,dV_g),
\ea
where $H^{-1}(U,\,dV_g)$ is the Sobolev space with smoothness $(-1)$
on $(U,g)$.
Since  $\gamma$ is a subset with (Hausdorff)
dimension 1
of the 3-dimensional domain $U$, it has zero capacitance.
Thus, the Lipschitz functions on $U$ that vanish on $\gamma$
are dense in $H^1(U)$, see \cite{KKM}.
Therefore, there  are no non-zero  distributions in $H^{-1}(U)$ supported on
$\gamma$. Hence
             we see that
\ba
\nabla\times E^e -ik \mu(x) H^e=0,\quad
\nabla\times H^e +ik \e(x) E^e=0\quad \hbox{in } U.
\ea
This also implies that
\ba
\nabla\cdotp (\e E^e)=0,\quad \nabla\cdotp (\mu H^e)=0\quad \hbox{in } U,
\ea
which, by  elliptic regularity, imply that $E^e$ and
$H^e$ are  $C^\infty$ smooth in $U$.

In summary, $E$ and $H$ have unique continuous extensions
to $\gamma$, and the extensions are classical
solutions to Maxwell's equations.
\hfill\proofbox

\bibliographystyle{amsalpha}

\vskip.2in

\noindent{\sc Department of Mathematics}

\noindent{\sc University of Rochester}

\noindent{\sc Rochester, NY 14627, USA, \quad \tt{allan@math.rochester.edu}}

\vskip.2in

\noindent{\sc Department of Mathematical Sciences}

\noindent{\sc University of Loughborough}

\noindent{\sc Loughborough, LE11 3TU, UK,\quad \tt{Y.V.Kurylev@lboro.ac.uk}}

\vskip.2in

\noindent{\sc Institute of Mathematics}

\noindent{\sc Helsinki University of Technology}

\noindent{\sc Espoo, FIN-02015, Finland,\quad \tt{Matti.Lassas@tkk.fi}}

\vskip.2in

\noindent{\sc Department of Mathematics}

\noindent{\sc University of Washington}

\noindent{\sc Seattle, WA 98195, USA,\quad \tt{gunther@math.washington.edu}}

\end{document}